\documentclass[a4paper,12pt]{article}
\usepackage{amssymb,fullpage,amsmath,amsthm,graphicx}
\newtheorem{thm}{Theorem}[section]
\newtheorem{lem}[thm]{Lemma}

\newtheorem{cor}[thm]{Corollary}

\begin{document}

\title{Random walk on the range of random walk}
\author{David A. Croydon}
\date{30 June 2009}
\maketitle

\begin{abstract}
We study the random walk $X$ on the range of a simple random walk on $\mathbb{Z}^d$ in dimensions $d\geq 4$. When $d\geq 5$ we establish quenched and annealed scaling limits for the process $X$, which show that the intersections of the original simple random walk path are essentially unimportant. For $d=4$ our results are less precise, but we are able to show that any scaling limit for $X$ will require logarithmic corrections to the polynomial scaling factors seen in higher dimensions. Furthermore, we demonstrate that when $d=4$ similar logarithmic corrections are necessary in describing the asymptotic behaviour of the return probability of $X$ to the origin.
\end{abstract}

\section{Introduction}

Let $S=(S_n)_{n\geq 0}$ be the simple random walk on $\mathbb{Z}^d$ starting from 0, built on an underlying probability space $\Omega$ with probability measure $\mathbf{P}$. Define the range of the random walk $S$ to be the graph $\mathcal{G}=(V(\mathcal{G}),E(\mathcal{G}))$ with vertex set
\[V(\mathcal{G}):=\left\{S_n:n\geq 0\right\},\]
and edge set
\[E(\mathcal{G}):=\left\{\{S_n,S_{n+1}\}:n\geq 0\right\}.\]
For $\mathbf{P}$-a.e. random walk path, the graph $\mathcal{G}$ is infinite, connected and clearly has bounded degree. In this article, the main object of study will be the discrete time simple random walk on $\mathcal{G}$, which we now introduce. For a given realisation of $\mathcal{G}$, write
\[X=\left((X_n)_{n\geq0},\mathbf{P}_x^\mathcal{G},x\in V(\mathcal{G})\right)\]
to represent the Markov chain with transition probabilities
\[P_{\mathcal{G}}(x,y):=\frac{1}{\mathrm{deg}_\mathcal{G}(x)}\mathbf{1}_{\{\{x,y\}\in E(\mathcal{G})\}},\hspace{20pt}\forall x,y\in V(\mathcal{G}),\]
where $\mathrm{deg}_\mathcal{G}(x)$ is the usual graph degree of $x$ in $\mathcal{G}$. For $x\in V(\mathcal{G})$, the law $\mathbf{P}_x^\mathcal{G}$ is the quenched law of the simple random walk on $\mathcal{G}$ started from $x$. Since $0$ is always an element of $V(\mathcal{G})$, we can also define an annealed law $\mathbb{P}$ for the random walk on $\mathcal{G}$ started from 0 as the semi-direct product of the environment law $\mathbf{P}$ and the quenched law $\mathbf{P}_0^\mathcal{G}$ by setting
\[\mathbb{P}:=\int \mathbf{P}_0^\mathcal{G}(\cdot){\rm d}\mathbf{P}.\]

When $d=1,2$ the recurrence of the random walk $S$ easily implies that $\mathcal{G}$ is $\mathbf{P}$-a.s. equal to the vertex set $\mathbb{Z}^d$ equipped with edges connecting points a unit Euclidean distance apart. Consequently, the law of $X$ under $\mathbb{P}$ and also under $\mathbf{P}_0^\mathcal{G}$, $\mathbf{P}$-a.s., is identical to the law of $S$ under $\mathbf{P}$, which is well-understood. In particular, it follows that $(n^{-1/2}X_{\lfloor nt\rfloor})_{t\geq 0}$ converges in distribution to standard Brownian motion in $\mathbb{R}^d$, and there are Gaussian bounds for the transition density of $X$. Conversely, for $d\geq 3$ the random walk $S$ is transient and does not explore all of $\mathbb{Z}^d$. In this case, since $\mathcal{G}$ has a non-deterministic structure, it becomes an interesting problem to determine the behaviour of $X$.

The problem of establishing dynamical properties of random walk paths has previously been investigated by physicists, with one motivation for doing so being its application to the study of the transport properties of sedimentary rocks of low porosity, where the commonly considered sub-critical percolation model does not reflect the pore connectivity properties seen in experiments, see \cite{BHK}, for example (references to other related models appear in \cite{HBA}, Section 8.4). In particular, numerical simulations have been conducted in an attempt to determine the walk dimension $d_W$ and spectral dimension $d_S$ of $X$, which are the exponents satisfying $\mathbb{E}|X_n|^2\approx n^{2/d_W}$, where $\mathbb{E}$ is the expectation under $\mathbb{P}$, and $\mathbb{P}(X_{2n}=0)\approx n^{-d_S/2}$ respectively \cite{HWBM} (see below for further discussion of these dimensions). Although we do not investigate here the case $d=3$, which was a main focus of \cite{HWBM}, our results do contribute to the existing higher-dimensional literature. More specifically, we will prove precise scaling results for $X$ and a quenched limit expression for the spectral dimension of $X$ when $d\geq5$, and also demonstrate that logarithmic corrections to these scaling results are necessary when $d=4$ .

To analyse the random walk on $\mathcal{G}$ when $d\geq 5$ it transpires that it is useful to introduce a second simple random walk $S'$ on $\mathbb{Z}^d$, which starts from 0 and is independent of $S$, and define a two-sided walk $\tilde{S}=(\tilde{S}_n)_{n\in\mathbb{Z}}$ by setting
\[\tilde{S}_n:=\left\{\begin{array}{ll}
                 S_n, & \mbox{if }n\geq 0, \\
                 S'_{-n}, &  \mbox{if }n<0.
               \end{array}\right.\]
It is known that the set of cut-times for the two-sided process $\tilde{S}$, which is defined by
\[\tilde{\mathcal{T}}:=\left\{n:\tilde{S}_{(-\infty,n]}\cap\tilde{S}_{[n+1,\infty)}=\emptyset\right\},\]
is infinite, $\mathbf{P}$-a.s., and, moreover, the point process of cut-times is stationary. This observation was applied in \cite{BSZ} to determine properties of a random walk in a particular high-dimensional random environment. In our case, by considering the sections of the random walk path between these cut-times, we obtain that $\tilde{\mathcal{G}}$, the graph defined from $\tilde{S}$ analogously to the definition of $\mathcal{G}$ from $S$, can be constructed by stringing together a stationary ergodic sequence of finite graphs (see Section \ref{ergsec}), and exploiting this decomposition of $\tilde{\mathcal{G}}$ we are able to determine the correct scaling for the random walk on $\tilde{\mathcal{G}}$ (see Theorem \ref{d5thmd}). Restricting to $\mathcal{G}$, we can subsequently deduce the following quenched and annealed scaling limits for the simple random walk $X$ on $\mathcal{G}$. The processes $B=(B_t)_{t\geq 0}$ and $W^{(d)}=(W^{(d)}_t)_{t\geq 0}$ are assumed to be independent standard Brownian motions on $\mathbb{R}$ and $\mathbb{R}^d$ respectively, both started from the origin. The notation $d_\mathcal{G}$ is used to represent the shortest path graph distance on $\mathcal{G}$.

{\thm \label{d5thm} Let $d\geq 5$. There exists a deterministic constant $\kappa_1{(d)}\in(0,\infty)$ such that, for $\mathbf{P}$-a.e. realisation of $\mathcal{G}$, the law of
\[\left(n^{-1/2}d_\mathcal{G}(0,X_{\lfloor tn\rfloor})\right)_{t\geq 0},\]
under $\mathbf{P}_0^{\mathcal{G}}$, converges as $n\rightarrow\infty$ to the law of $(|B_{t\kappa_1{(d)}}|)_{t\geq 0}$. There also exists a deterministic constant $\kappa_2{(d)}\in(0,\infty)$ such that the law of
\[\left(n^{-1/4}X_{\lfloor tn\rfloor}\right)_{t\geq 0},\]
under $\mathbb{P}$, converges as $n\rightarrow\infty$ to the law of $(W^{(d)}_{|B_{t\kappa_2{(d)}}|})_{t\geq 0}$.}
\bigskip

By imitating the construction in \cite{Croydonhaus} of the Brownian motion on the range of the super-Brownian motion conditioned to have total mass equal to one, it is possible to interpret the limiting process of part (b) of the above theorem as the Brownian motion on the range $\mathcal{R}:=\{W_t^{(d)}:t\geq 0\}$ of the Brownian motion $W^{(d)}$. In particular, if $d\geq 4$, $\mathbf{P}$-a.s., the map $t\mapsto W_t^{(d)}$ is injective, and so between any two points $x,y\in\mathcal{R}$, there is a unique arc. Moreover, by applying \cite{Cies}, Theorem 5, $\mathbf{P}$-a.s., we can define a metric $d_\mathcal{R}$ on $\mathcal{R}$ by setting $d_\mathcal{R}(x,y)$ to be equal to the Hausdorff measure with gauge function $c(d) x^2 \ln\ln x^{-1}$ of the arc between $x$ and $y$ in $\mathcal{R}$, where the deterministic constant $c(d)$, depending only on $d$, can be chosen so that if $x=W^{(d)}_s$ and $y=W_t^{(d)}$, then $d_\mathcal{R}(x,y)=|s-t|$. Since $t\mapsto W^{(d)}_t$ is a measure preserving isometry from $\mathbb{R}_+$, equipped with the one-dimensional Lebesgue measure, to $(\mathcal{R},d_\mathcal{R},\mu^\mathcal{R})$, where $\mu^\mathcal{R}$ is the Hausdorff measure with gauge function $c(d) x^2 \ln\ln x^{-1}$ restricted to $\mathcal{R}$, there is no problem in checking that the canonical Brownian motion $B^\mathcal{R}$ on the metric-measure space $(\mathcal{R},d_\mathcal{R},\mu^\mathcal{R})$ is simply the image under the map $t\mapsto W^{(d)}_t$ of the standard Brownian motion on $\mathbb{R}_+$ reflected at 0, for $\mathbf{P}$-a.e. realisation of $\mathcal{R}$ (see \cite{Aldous2}, Section 5.2, for a natural definition of a Brownian motion on a tree-like metric-measure space). After checking that the law of $B^\mathcal{R}$ can constructed in a $W^{(d)}$-measurable way, which is much more straightforward than the corresponding measurability result of \cite{Croydonhaus}, Proposition 7.2, in dimensions $d\geq 4$ we are able to define the annealed law of $B^\mathcal{R}$ by averaging out over the law of $W^{(d)}$. The resulting $C(\mathbb{R}_+,\mathbb{R}^d)$-valued process, the Brownian motion on the range of Brownian motion, satisfies
\[\left(B^\mathcal{R}_t\right)_{t\geq 0}=\left(W^{(d)}_{|B_{t}|}\right)_{t\geq 0}\]
in distribution, and is therefore, when $d\geq 5$, the scaling limit of the simple random walk on the range of a simple random walk, at least up to a constant time-change.

For $d\geq 5$ further understanding of the random walk on the range of random walk is provided by the following result, which demonstrates that the quenched probability that the process $X$ returns to the origin after $2n$ steps decays polynomially with the same exponent as for standard simple random walk on $\mathbb{Z}$. As with the previous result, its proof relies on determining properties of the one-sided graph $\mathcal{G}$ using the ergodic description we have for the two-sided graph $\tilde{\mathcal{G}}$. See (\ref{anntp}) for the corresponding annealed bounds.

{\thm \label{d5hk}Let $d\geq 5$. There exist deterministic constants $c_1,c_2\in(0,\infty)$ such that, for $\mathbf{P}$-a.e. realisation of $\mathcal{G}$,
\begin{equation}\label{qtp}
c_1n^{-1/2}\leq \mathbf{P}^\mathcal{G}_0(X_{2n}=0)\leq c_2n^{-1/2},
\end{equation}
for large $n$.}
\bigskip

For $d=4$ the set of cut-times $\tilde{\mathcal{T}}$ is empty, $\mathbf{P}$-a.s., and so we can not use the same ergodic arguments to analyse $X$ as in the higher-dimensional case. However, the one-sided simple random walk $S$ still admits cut-times, by which we mean that the set
\begin{equation}\label{osct}
\mathcal{T}:=\left\{n:\tilde{S}_{[0,n]}\cap\tilde{S}_{[n+1,\infty)}=\emptyset\right\},
\end{equation}
is non-empty, $\mathbf{P}$-a.s., and considering the structure of $\mathcal{G}$ between the cut-points $(S_n)_{n\in\mathcal{T}}$ will prove similarly helpful in understanding the random walk $X$ on $\mathcal{G}$. Whilst we are unable to prove exact scaling in this dimension, we are able to show that any such result will require extra logarithmic correction terms compared to higher dimensions, thereby demonstrating that the random walk is anomalous in this dimension (see \cite{HBA} for an excellent survey of work regarding anomalous diffusions in disordered media). Thus, in the language of statistical mechanics, this result establishes that the critical dimension of the random walk on the range of the random walk is 4. We note that the fact that the intersections of the of the original simple random walk affect the behaviour in a logarithmic way in $d=4$, but can effectively be neglected for higher dimensions, is not a surprise given the analogous results known to hold for other random walk models depending on self-interaction properties. These include the self-avoiding walk, for which logarithmic corrections have not yet been rigourously proved in $4$ dimensions, see \cite{Slade}, Chapter 2, for a summary of mathematical results, but have been observed in physics, where relevant work includes \cite{Brezin}, \cite{deg}, \cite{Dup}, \cite{Grassberger}, and also the loop-erased random walk introduced by Lawler \cite{LawO} (see also \cite{log}). In fact, in proving the subsequent theorem, we exploit the fact that the number of steps in the loop-erasure of the path from $S_m$ to $S_n$ gives an upper bound for graph distance in $\mathcal{G}$ between $S_m$ and $S_n$.

{\thm\label{d4thm} If $d=4$, then
\begin{equation}\label{d4lim}
\lim_{\lambda\rightarrow\infty}\liminf_{n\rightarrow\infty}\mathbb{P}\left(\lambda^{-1} n^{1/4}(\ln n)^{1/24}\leq \max_{m\leq n}|X_m|\leq \lambda n^{1/4} (\ln n)^{7/12}\right)=1.
\end{equation}}

In addition to the above result regarding the scaling of the random walk $X$ on $\mathcal{G}$, when $d=4$ we can also establish bounds for the quenched transition density of the random walk on $\mathcal{G}$, which confirms that logarithmic corrections to the bounds that hold in higher dimensions are necessary for this quantity too.

{\thm \label{tdthm} If $d=4$, then
\[\lim_{\lambda\rightarrow\infty}\liminf_{n\rightarrow\infty}\mathbb{P}\left(\lambda^{-1} n^{-1/2}(\ln n)^{-3/2}\leq \mathbf{P}^\mathcal{G}_0(X_{2n}=0)\leq \lambda n^{-1/2} (\ln n)^{-1/6}\right)=1.\]}

Let us continue by observing that the results of Theorems \ref{d5thm} and \ref{d4thm} (and also Lemma \ref{d4exit} below) imply a distributional version of the result that that the walk dimension $d_W$ of $X$ is 4 with respect to the Euclidean distance and $2$ with respect to the graph distance when $d\geq 4$, which contrasts with the low-dimensional $d\leq 2$ setting, where $d_W=2$ with respect to both distances. Furthermore, defining the quenched spectral dimension of the random walk $X$ by the limit
\[d_S:=\lim_{n\rightarrow\infty}\frac{2\ln\mathbf{P}^\mathcal{G}_0(X_{2n}=0)}{-\ln n},\]
when it exists, our results demonstrate that the spectral dimension depends on $d$ in the following way:  $d_S= 1$, $\mathbf{P}$-a.s., for $d=1$ and $d\geq 5$;  $d_S= 2$, $\mathbf{P}$-a.s., for $d=2$; $d_S$ defined by a probabilistic limit is equal to 1 for $d=4$ (in a recent preprint, this has actually been shown to be true $\mathbf{P}$-a.s. \cite{Shira}). The problem of determining $d_W$ and $d_S$ when $d=3$ seems difficult, and we do not present any progress on this problem here, but merely remark that physicists' numerical simulations suggest that $d_W\approx 7/2$ with respect to the Euclidean distance, and $d_S\approx 8/7$ in this dimension \cite{HWBM}, and note that a non-trivial bound confirming that $d_S>1$ is proved in \cite{Shira}.

Finally, another natural choice of transition probabilities for the random walk $X$ on $\mathcal{G}$ is to set $P_\mathcal{G}(x,y)=\mu_{xy}/\mu_x$ for $\{x,y\}\in E(\mathcal{G})$, where
\begin{equation}\label{muxy}
\mu_{xy}:=\#\left\{n:\{S_n,S_{n+1}\}=\{x,y\}\right\}
\end{equation}
is the number of crossings of the edge $\{x,y\}$ by $S$ and $\mu_x:=\sum_{y:\{x,y\}\in E(\mathcal{G})}\mu_{xy}$, so that the random walk $X$ is more likely to jump along edges that the random walk $S$ traversed more frequently. With straightforward modifications to the proofs, all the results stated in the introduction, albeit with suitably adjusted constants, will hold for this random walk. See the end of Section \ref{slsec} for elaboration on this point.

The article is organised as follows. In Section \ref{ergsec}, we describe the decomposition of $\tilde{\mathcal{G}}$ into a stationary ergodic sequence of finite graphs for $d\geq 5$, and then proceed in Section \ref{slsec} to apply this to proving Theorems \ref{d5thm} and \ref{d5hk}. Finally, in Section \ref{d4sec}, we investigate the behaviour of $X$ in the critical dimension $d=4$ in order to establish Theorems \ref{d4thm} and \ref{tdthm}.

\section{Ergodic behaviour for $d\geq 5$}\label{ergsec}

In this section we suppose that $d\geq 5$ and study the two-sided graph $\tilde{\mathcal{G}}$. Let us start by making explicit our probability space by setting $\Omega$ to be equal to the countable product of discrete spaces $\{x\in\mathbb{Z}^d:|x|=1\}^\mathbb{Z}$, and equipping this set with the product $\sigma$-algebra. We suppose that $\mathbf{P}$ is defined to satisfy
\[\mathbf{P}((\omega_{-m},\dots,\omega_n)=(x_{-m},\dots,x_n))=(2d)^{-(m+n+1)},\]
for every $(x_{-m},\dots,x_n)\in \{x\in\mathbb{Z}^d:|x|=1\}^{m+n+1}$, $m,n\geq 0$, and also that
\[\tilde{S}_n=\tilde{S}_n(\omega):=\left\{\begin{array}{ll}
                                            \omega_1+\dots\omega_n, & n\geq 1, \\
                                            0, & n=0, \\
                                            -(\omega_{n+1}+\dots+\omega_0), & n\leq -1.
                                          \end{array}\right.\]
We will denote by $(\theta_n)_{n\in\mathbb{Z}}$ the canonical shift maps on $\Omega$, so that $(\tilde{S}\circ\theta_m)_n=\tilde{S}_{m+n}-\tilde{S}_m$ for $m,n\in\mathbb{Z}$.

In the high-dimensional case we are considering here, it is possible to check that there is a strictly positive probability that $0$ is a cut-time for $\tilde{S}$ by applying results of \cite{ET}, and consequently the probability measure $\hat{\mathbf{P}}:=\mathbf{P}(\cdot|0\in\tilde{\mathcal{T}})$ is well-defined. Furthermore, also by \cite{ET}, the measurable set $\Omega^*\subseteq\Omega$ where $\tilde{\mathcal{T}}\cap(-\infty,0]$ and $\tilde{\mathcal{T}}\cap[0,\infty)$ are both infinite satisfies $\mathbf{P}(\Omega^*)=1$, and we henceforth suppose that $\mathbf{P}$ and $\hat{\mathbf{P}}$ are restricted to this set. In particular, we can always write $\tilde{\mathcal{T}}=\{T_n:n\in\mathbb{Z}\}$, where $\dots T_{-2}<T_{-1}<T_0\leq0<T_1<T_2<\dots$. Fundamental in proving many of our subsequent results is the following lemma, where $\mathbf{E}$ and $\hat{\mathbf{E}}$ are the expectations under $\mathbf{P}$ and $\hat{\mathbf{P}}$ respectively.

{\lem The measure $\hat{\mathbf{P}}$ is invariant under and ergodic for $\hat{\theta}:=\theta_{T_1}$. Moreover, for bounded measurable $f$,
\begin{equation}\label{pchar}
{\mathbf{E}}(f)=\frac{\hat{\mathbf{E}}\left(\sum_{n=0}^{T_1-1}f\circ\theta_n \right)}{\hat{\mathbf{E}}(T_1)},
\end{equation}
where $\hat{\mathbf{E}}(T_1)\in[1,\infty)$.}
\begin{proof} See \cite{BSZ}, Lemma 1.1 and proof of \cite{BSZ}, Proposition 1.3.
\end{proof}

If we define for each $n\in\mathbb{Z}$ a graph $\tilde{\mathcal{G}}_n$ to have vertex set
\[V(\tilde{\mathcal{G}}_n):=\left\{\tilde{S}_m-C_n:T_n\leq m\leq T_{n+1}\right\},\]
and edge set
\[E(\tilde{\mathcal{G}}_n):=\left\{\{\tilde{S}_{m}-C_n,\tilde{S}_{m+1}-C_n\}:T_n\leq m< T_{n+1}\right\},\]
where we write $C_n$ to represent the cut-point $\tilde{S}_{T_n}$, then the previous lemma immediately implies that the sequence of finite graphs with distinguished vertices $((\tilde{\mathcal{G}}_n, C_{n+1}-C_n))_{n\in\mathbb{Z}}$ is stationary and ergodic under $\hat{\mathbf{P}}$. Moreover, it is clear that $\tilde{\mathcal{G}}$ can be reconstructed from $((\tilde{\mathcal{G}}_n, C_{n+1}-C_n))_{n\in\mathbb{Z}}$ by adjoining graphs at cut-points. As a consequence of this ergodicity, we are able to define finitely a number of $\hat{\mathbf{P}}$-a.s. limits and expectations under $\hat{\mathbf{P}}$. The quantities defined in the following lemma will later appear in the definitions of the diffusion constants for the scaling limits of the random walks on $\tilde{\mathcal{G}}$ and $\mathcal{G}$. We write $d_{\tilde{\mathcal{G}}}$ to represent the shortest path graph distance on $\tilde{\mathcal{G}}$. The function $R_{\tilde{\mathcal{G}}}$ is the effective resistance on $\tilde{\mathcal{G}}$ when we suppose that a unit resistor is placed along each edge (see \cite{Barlow}, Definition 4.23, for example). The usual graph degree of $x$ in $\tilde{\mathcal{G}}$ is denoted by $\mathrm{deg}_{\tilde{\mathcal{G}}}(x)$.

{\lem \label{exponent}$\hat{\mathbf{P}}$-a.s., we have that
\begin{equation}\label{taud}
\frac{T_n}{n}\rightarrow \tau(d):=\hat{\mathbf{E}}T_1\in[1,\infty),
\end{equation}
\begin{equation}\label{deltad}
\frac{d_{\tilde{\mathcal{G}}}(0,C_n)}{|n|}\rightarrow \delta(d):=\hat{\mathbf{E}}d_{\tilde{\mathcal{G}}}(0,C_1)\in[1,\infty),
\end{equation}
\begin{equation}\label{rhod}
\frac{R_{\tilde{\mathcal{G}}}(0,C_n)}{|n|}\rightarrow\rho(d):=\hat{\mathbf{E}}R_{\tilde{\mathcal{G}}}(0,C_1)\in[1,\infty),
\end{equation}
as $|n|\rightarrow\infty$. Furthermore,
\[\nu(d):=\tfrac{1}{2}\hat{\mathbf{E}}{\rm{deg}}_{\tilde{\mathcal{G}}} (0)\in[1,\infty).\]}
\begin{proof} If $n\in\mathbb{N}$, then $T_n-T_0=\sum_{m=1}^n (T_{m}-T_{m-1})$. Thus the limit at (\ref{taud}) as $n\rightarrow \infty$ results from the ergodic theorem and the finiteness of $\tau(d)$, which was noted in the previous lemma. Similarly for $n\rightarrow -\infty$. To establish (\ref{deltad}), observe that $|n|\leq d_{\tilde{\mathcal{G}}}(0,C_n)\leq |T_n|$ for every $n\in\mathbb{Z}$, from which it follows that $\limsup_{|n|\rightarrow\infty}d_{\tilde{\mathcal{G}}}(0,C_n)/|n|\in[1,\tau(d)]$, $\hat{\mathbf{P}}$-a.s. Thus, the ergodic theorem again implies the limiting result with $\delta(d)\in[1,\tau(d)]$. The resistance on a graph is always bounded above by the graph distance. Hence the expectation in (\ref{rhod}) is finite and satisfies $\rho(d)\leq\delta(d)$. Moreover, whenever $0\in\tilde{\mathcal{T}}$, it is clear that any path in $\tilde{\mathcal{G}}$ from $0$ to $C_1$ must contain the edge $\{S_0,S_1\}$. This implies that $R_{\tilde{\mathcal{G}}}(0,C_1)\geq 1$, $\hat{\mathbf{P}}$-a.s., and therefore $\rho(d)\geq 1$.  The series law for resistors allows us to again apply the ergodic theorem to deduce the limit exists. Finally, it is elementary to check that the remaining expectation satisfies $1\leq\nu(d)\leq\tau(d)$, which completes the proof.
\end{proof}

\section{Scaling limit and transition probability for $d\geq 5$}\label{slsec}

Applying the description of $\tilde{\mathcal{G}}$ from the previous section, we now proceed to analyse the associated random walk, with the first aim of this section being to prove a two-sided version of Theorem \ref{d5thm}. We will write $\tilde{X}$ to represent the random walk on $\tilde{\mathcal{G}}$; its quenched law started from 0 will be denoted $\mathbf{P}^{\tilde{\mathcal{G}}}_0$. Define the hitting times by $\tilde{X}$ of the set of cut-points $\tilde{\mathcal{C}}:=\{C_n:n\in\mathbb{Z}\}$ by
\begin{equation}\label{hdef1}
H_0:=\min\{m\geq 0:\tilde{X}_m\in\tilde{\mathcal{C}}\},
\end{equation}
and, for $n\geq 1$,
\begin{equation}\label{hdef2}
H_n:=\min\{m> H_{n-1}:\tilde{X}_m\in\tilde{\mathcal{C}}\}.
\end{equation}
Denote by $\pi$ the bijection from $\mathbb{Z}$ to $\tilde{\mathcal{C}}$ that satisfies $\pi(n)=C_n$, and let $J=(J_n)_{n\geq 0}$ be the $\mathbb{Z}$-valued process obtained by setting
\[J_n:=\pi^{-1}\left(\tilde{X}_{H_n}\right).\]
Note that $J$ can remain at a particular integer for multiple time-steps. Using techniques developed for the random conductance model, it is possible to deduce the following quenched scaling limit for $J$.

{\lem \label{jump}  Let $d\geq 5$. For $\hat{\mathbf{P}}$-a.e. realisation of $\tilde{\mathcal{G}}$, the law of $(n^{-1/2}J_{\lfloor tn\rfloor})_{t\geq 0}$ under $\mathbf{P}_0^{\tilde{\mathcal{G}}}$ converges to the law of $(B_{t/\nu(d)\rho(d)})_{t\geq 0}$.}
\begin{proof} Fix a particular $\omega\in\Omega^*\cap\{0\in\tilde{\mathcal{T}}\}$, set $\tilde{\mathcal{G}}=\tilde{\mathcal{G}}(\omega)$ and let $n\in\mathbb{Z}$. Conditional on $J_m=n$, the amount of time that $J$ spends at position $n$ from time $m$ until it next hits $\{n-1,n+1\}$ is geometric, parameter $p$, where $p:=\mathbf{P}^{\tilde{\mathcal{G}}}_0(J_{m+1}=n|J_m=n)$, and therefore has expectation $(1-p)^{-1}$. Applying the definition of $J$, it is possible to check that this expectation is also equal to
\begin{eqnarray*}
\mathbf{E}^{\tilde{\mathcal{G}}}_{C_n}\left(\sum_{m=0}^{H(\{C_{n-1},C_{n+1}\})-1}\mathbf{1}_{\{\tilde{X}_m=C_n\}}\right)&=&{\rm deg}_{\tilde{\mathcal{G}}}(C_n)R_{\tilde{\mathcal{G}}}(C_n,\{C_{n-1},C_{n+1}\})\\
&=&\frac{{\rm deg}_{\tilde{\mathcal{G}}}(C_n) R_{\tilde{\mathcal{G}}}(C_{n-1},C_n)R_{\tilde{\mathcal{G}}}(C_n,C_{n+1})}{R_{\tilde{\mathcal{G}}}(C_{n-1},C_{n+1})},
\end{eqnarray*}
where $H(\{C_{n-1},C_{n+1}\}):=\min\{m:\tilde{X}_m\in\{C_{n-1},C_{n+1}\}\}$, the first equality is an application of a well-known electrical network interpretation of the occupation density of a killed random walk on a graph  (see \cite{DS} or \cite{LP}, Chapter 2, for example)
and the second equality follows from the parallel law for electrical resistance. Rearranging for $p$, we find that
\begin{equation}\label{same}
\mathbf{P}^{\tilde{\mathcal{G}}}_0(J_{m+1}=n|J_m=n)=1-\frac{R_{\tilde{\mathcal{G}}}(C_{n-1},C_{n+1})}{{\rm deg}_{\tilde{\mathcal{G}}}(C_n) R_{\tilde{\mathcal{G}}}(C_{n-1},C_n)R_{\tilde{\mathcal{G}}}(C_n,C_{n+1})}.
\end{equation}
Further elementary calculations allow it to be deduced that
\begin{equation}\label{notsame}
\mathbf{P}^{\tilde{\mathcal{G}}}_0(J_{m+1}=n\pm1|J_m=n)=\frac{1}{{\rm deg}_{\tilde{\mathcal{G}}}(C_n) R_{\tilde{\mathcal{G}}}(C_n,C_{n\pm1})}.
\end{equation}
These formulae easily imply that the process $({\rm sgn}(J_n)R_{\tilde{\mathcal{G}}}(0,\pi(J_n)))_{n\geq 0}$ is a martingale, and the stationary ergodic decomposition of $\tilde{\mathcal{G}}$ of the previous section allows us to apply the Lindeberg-Feller central limit theorem to this martingale by making only simple adaptations to the ``environment viewed from the particle'' argument described in the introduction of \cite{BisP} for a random walk among stationary ergodic random conductances. In particular, to deduce that the law of $(n^{-1/2}{\rm sgn}(J_{\lfloor tn\rfloor})R_{\tilde{\mathcal{G}}}(0,\pi(J_{\lfloor tn\rfloor})))_{t\geq 0}$ under $\mathbf{P}^{\tilde{\mathcal{G}}}_0$ converges to a Brownian motion law as $n\rightarrow\infty$, it will suffice to demonstrate the square integrability condition
\[\frac{\hat{\mathbf{E}}\left({\rm deg}_{\tilde{\mathcal{G}}}(0)\mathbf{E}_0^{\tilde{\mathcal{G}}}R_{\tilde{\mathcal{G}}}(0,\tilde{X}_{H_1})^2\right)}{
\hat{\mathbf{E}}{\rm deg}_{\tilde{\mathcal{G}}}(0)}<\infty\]
holds. Moreover, we will show that the left-hand side above, which provides the limiting diffusion constant, is equal to $\rho(d)/\nu(d)$. Note that the ${\rm deg}_{\tilde{\mathcal{G}}}(0)/\hat{\mathbf{E}}{\rm deg}_{\tilde{\mathcal{G}}}(0)$ factor arises here as a result of the fact that the invariant measure of $J$ is given by $({\rm deg}_{\tilde{\mathcal{G}}}(C_n))_{n\in\mathbb{Z}}$, which can be checked using (\ref{notsame}) and the detailed-balance equations. From the transition probabilities at (\ref{same}) and (\ref{notsame}), we have that
\[{\rm deg}_{\tilde{\mathcal{G}}}(0)\mathbf{E}_0^{\tilde{\mathcal{G}}}R_{\tilde{\mathcal{G}}}(0,\tilde{X}_{H_1})^2=R_{\tilde{\mathcal{G}}}(0,C_1)+R_{\tilde{\mathcal{G}}}(0,C_{-1}).\]
Hence, by stationarity,
\[\hat{\mathbf{E}}\left({\rm deg}_{\tilde{\mathcal{G}}}(0)\mathbf{E}_0^{\tilde{\mathcal{G}}}R_{\tilde{\mathcal{G}}}(0,\tilde{X}_{H_1})^2\right)=2\rho(d),\]
which confirms that the limiting diffusion constant is indeed $\rho(d)/\nu(d)$. To complete the proof, we note from (\ref{rhod}) that ${\rm sgn}(n)R_{\tilde{\mathcal{G}}}(0,\pi(n))\sim n\rho(d)$ as $|n|\rightarrow\infty$, and hence the law of $(n^{-1/2}J_{\lfloor tn\rfloor})_{t\geq 0}$ under $\mathbf{P}_0^{\tilde{\mathcal{G}}}$ converges to the law of a Brownian motion with diffusion constant $\rho(d)/\nu(d)\rho(d)^2=1/\nu(d)\rho(d)$.
\end{proof}

We now show that the hitting times $(H_n)_{n\geq 0}$ defined at (\ref{hdef1}) and (\ref{hdef2}) grow linearly asymptotically. In the proof of this result, we consider the measure $\mu_{\tilde{\mathcal{G}}}$ on $V(\tilde{\mathcal{G}})$ defined to satisfy
\[\mu_{\tilde{\mathcal{G}}}(\{x\}):={{\rm deg}_{\tilde{\mathcal{G}}}(x)},\hspace{20pt}\forall x\in V(\tilde{G}).\]
It arises naturally in the argument, because it is invariant for the simple random walk $\tilde{X}$ on $\tilde{\mathcal{G}}$.

{\lem\label{hitting}  Let $d\geq 5$. For $\hat{\mathbf{P}}$-a.e. realisation of $\tilde{\mathcal{G}}$, $\mathbf{P}^{\tilde{\mathcal{G}}}_0$-a.s., we have that
\[\frac{H_n}{n}\rightarrow \eta(d):= \frac{\hat{\mathbf{E}}\left({\rm deg}_{\tilde{\mathcal{G}}}(0)\mathbf{E}^{\tilde{\mathcal{G}}}_0 H_1\right)}{\hat{\mathbf{E}}{\rm deg}_{\tilde{\mathcal{G}}}(0)}\in[1,\infty).\]}
\begin{proof}  We start by checking that $\eta(d)\in[1,\infty)$. Fix $\omega\in\Omega^*\cap\{0\in\tilde{\mathcal{T}}\}$ and set $\tilde{\mathcal{G}}=\tilde{\mathcal{G}}(\omega)$. Applying the Markov property of $\tilde{X}$ and standard bounds for hitting times of random walks on graphs in terms of resistance and volume (see \cite{Barlow}, Corollary 4.28, for example), we find
\begin{eqnarray}
1\leq \mathbf{E}_0^{\tilde{\mathcal{G}}}H_1&\leq &1+\sup_{\substack{x\in V(\tilde{\mathcal{G}})\backslash \{C_{-1},C_1\}:\\\{0,x\}\in E(\tilde{\mathcal{G}})}}\mathbf{E}_x^{\tilde{\mathcal{G}}}H_1\nonumber\\
&\leq &1+\sup_{\substack{x\in V(\tilde{\mathcal{G}})\backslash \{C_{-1},C_1\}:\\\{0,x\}\in E(\tilde{\mathcal{G}})}}R_{\tilde{\mathcal{G}}}(x,\{C_{-1},0,C_1\})\mu_{\tilde{\mathcal{G}}}\left(A\backslash\{0\}\right),\label{bound1}
\end{eqnarray}
where $A$ is the graph connected component of $V(\tilde{\mathcal{G}})\backslash\{C_{-1},C_1\}$ and we set the supremum of an empty set to be zero. Since any $x$ in the set over which above supremum is taken is connected to $\{C_{-1},0,C_1\}$ by a path consisting of a single edge, the resistance in this expression is bounded above by 1. A simple counting exercise also implies that
\begin{equation}\label{bound2}
\mu_{\tilde{\mathcal{G}}}\left(A\backslash\{0\}\right)\leq 2(T_{1}-T_{-1})-{\rm deg}_{\tilde{\mathcal{G}}}(0)-2\leq 2(T_1-T_{-1})-4.
\end{equation}
Therefore, taking expectations, $1\leq \eta(d)\leq 1+4d\tau(d)/\nu(d)<\infty$.

The limit statement of the lemma can be proved similarly to \cite{DeM}, (4.16), replacing the exponential holding times of that article with the random hitting times $(H_{n+1}-H_n)_{n\geq0}$. In particular, it is possible to check that the sequence \[((\tilde{\mathcal{G}}-\tilde{X}_{H_n},H_{n+1}-H_n))_{n\geq0}\]
is ergodic under the annealed measure $\int \mathbf{P}_0^{\tilde{\mathcal{G}}}(\cdot){\rm d}\hat{\mathbf{P}}$, where $\tilde{\mathcal{G}}-\tilde{X}_{H_n}$ is the graph with vertex set $\{x-\tilde{X}_{H_n}:x\in V(\tilde{\mathcal{G}})\}$ and edge set $\{\{x-\tilde{X}_{H_n},y-\tilde{X}_{H_n}\}:\{x,y\}\in E(\tilde{\mathcal{G}})\}$.
Moreover, from the expressions for the transition probabilities of $J$ at (\ref{same}) and (\ref{notsame}), it is possible to deduce that the invariant measure of $((\tilde{\mathcal{G}}-\tilde{X}_{H_n},H_{n+1}-H_n))_{n\geq0}$ is given by the law of $(\tilde{\mathcal{G}},H_1)$ under the size-biased measure $\int {\rm deg}_{\tilde{\mathcal{G}}}(0)\mathbf{P}_0^{\tilde{\mathcal{G}}}(\cdot){\rm d}\hat{\mathbf{P}}/\hat{\mathbf{E}}{\rm deg}_{\tilde{\mathcal{G}}}(0)$, which is clearly bi-absolutely continuous with respect to $\int \mathbf{P}_0^{\tilde{\mathcal{G}}}(\cdot){\rm d}\hat{\mathbf{P}}$. Consequently, the proof of the lemma can be completed by an application of the ergodic theorem.
\end{proof}

The subsequent corollary demonstrates that the two previous results are unaffected by dropping the conditioning on the event $\{0\in\tilde{\mathcal{T}}\}$.

{\cor\label{coro}  Let $d\geq 5$. Lemmas \ref{jump} and \ref{hitting} also hold for $\mathbf{P}$-a.e. realisation of $\tilde{\mathcal{G}}$.}
\begin{proof} Fix $\omega\in\Omega^*$ and set $\tilde{\mathcal{G}}=\tilde{\mathcal{G}}(w)$. Suppose that the conclusions of Lemmas \ref{jump} and \ref{hitting} hold for both $\tilde{\mathcal{G}}\circ\theta_{T_0}$ and $\tilde{\mathcal{G}}\circ\theta_{T_1}$. Noting that $H_0<\infty$, $\mathbf{P}_0^{\tilde{\mathcal{G}}}$-a.s., we can write
\[\mathbf{P}_0^{\tilde{\mathcal{G}}}=\mathbf{P}_0^{\tilde{\mathcal{G}}}(\cdot|J_0=0)\mathbf{P}_0^{\tilde{\mathcal{G}}}(J_0=0)+\mathbf{P}_0^{\tilde{\mathcal{G}}}(\cdot|J_0=1)\mathbf{P}_0^{\tilde{\mathcal{G}}}(J_0=1).\]
The law of $(J_n,H_n-H_0)_{n\geq 0}$ under $\mathbf{P}_0^{\tilde{\mathcal{G}}}(\cdot|J_0=0)$ is equal to the law of $(J_n,H_n)_{n\geq 0}$ under $\mathbf{P}_0^{\tilde{\mathcal{G}}\circ \theta_{T_0}}$. Similarly, the law of $(J_n,H_n-H_0)_{n\geq 0}$ under $\mathbf{P}_0^{\tilde{\mathcal{G}}}(\cdot|J_0=1)$ is equal to the law of $(J_n+1,H_n)_{n\geq 0}$ under $\mathbf{P}_0^{\tilde{\mathcal{G}}\circ \theta_{T_1}}$. It follows that the law of $(n^{-1/2}J_{\lfloor tn\rfloor })_{t\geq 0}$ under $\mathbf{P}_0^{\tilde{\mathcal{G}}}$ converges to that of a Brownian motion with diffusion constant $1/\nu(d)\rho(d)$, and $H_n/n$ converges to $\eta(d)$, $\mathbf{P}_0^{\tilde{\mathcal{G}}}$-a.s. To complete the proof, note that the characterisation of $\mathbf{P}$ at (\ref{pchar}), combined with the invariance of $\hat{\mathbf{P}}$ under $\hat{\theta}=\theta_{T_1}$, allows us to apply Lemmas \ref{jump} and \ref{hitting} to deduce that our assumptions on $\tilde{\mathcal{G}}$ hold $\mathbf{P}$-a.s.
\end{proof}

This result allows us to prove our two-sided version of Theorem \ref{d5thm} with
\[\kappa_1(d):=\frac{{\delta(d)}^2}{\nu(d)\rho(d)\eta(d)}\in(0,\infty),\]
\[\kappa_2(d):=\frac{{\tau(d)}^2}{\nu(d)\rho(d)\eta(d)}\in(0,\infty).\]
We define the annealed measure of the random walk $\tilde{X}$ on $\tilde{\mathcal{G}}$ by setting $\tilde{\mathbb{P}}:=\int\mathbf{P}_0^{\tilde{\mathcal{G}}}(\cdot){\rm d}\mathbf{P}$.

{\thm \label{d5thmd} Let $d\geq 5$. For $\mathbf{P}$-a.e. realisation of $\tilde{\mathcal{G}}$, the law of
\[\left(n^{-1/2}d_{\tilde{\mathcal{G}}}(0,\tilde{X}_{\lfloor tn\rfloor})\right)_{t\geq 0},\]
under $\mathbf{P}_0^{\tilde{\mathcal{G}}}$, converges as $n\rightarrow\infty$ to the law of $(|B_{t\kappa_1{(d)}}|)_{t\geq 0}$. Furthermore, the law of
\[\left(n^{-1/4}\tilde{X}_{\lfloor tn\rfloor}\right)_{t\geq 0},\]
under $\tilde{\mathbb{P}}$, converges as $n\rightarrow\infty$ to the law of $(W^{(d)}_{B_{t\kappa_2{(d)}}})_{t\geq 0}$.}
\begin{proof} If the discrete time inverse $H^{-1}$ of $H$ is defined by setting
\[H^{-1}_n:=\min\{m:H_m>n\},\]
then Corollary \ref{coro} implies that the law of \[\left(n^{-1/2}\delta(d)\left|Z_{\lfloor tn\rfloor}\right|\right)_{t\geq0}\]
under $\mathbf{P}_0^{\tilde{\mathcal{G}}}$ converges to the law  of $(|B_{t\kappa_1{(d)}}|)_{t\geq 0}$, for $\mathbf{P}$-a.e. realisation of $\tilde{\mathcal{G}}$, where $Z=(Z_n)_{n\geq 0}$ is defined by setting
\begin{equation}\label{zdef}
Z_n:=J_{H^{-1}_{n}}.
\end{equation}
Hence, we will obtain that the same convergence result holds for $(n^{-1/2}d_{\tilde{\mathcal{G}}}(0,\tilde{X}_{\lfloor tn\rfloor}))_{t\geq 0}$, if we can show that, for every $\varepsilon,T\in(0,\infty)$,
\begin{equation}\label{close}
\lim_{n\rightarrow\infty}\mathbf{P}_0^{\tilde{\mathcal{G}}}\left(n^{-1/2}\sup_{t\in[0,T]} \left|d_{\tilde{\mathcal{G}}}(0,\tilde{X}_{\lfloor tn\rfloor})-\delta(d)|Z_{\lfloor tn\rfloor}|\right|>\varepsilon\right)=0,
\end{equation}
for $\mathbf{P}$-a.e. $\tilde{\mathcal{G}}$. Writing $Z_n^*:=\sup_{m\leq {\lfloor Tn\rfloor}}|Z_m|$, the definitions of $\pi$ and $Z$ imply that
\begin{eqnarray*}
\lefteqn{\sup_{t\in[0,T]}\left|d_{\tilde{\mathcal{G}}}(0,\tilde{X}_{\lfloor tn\rfloor})-\delta(d)|Z_{\lfloor tn\rfloor}|\right|}\\
&\leq&
\sup_{t\in[0,T]}\left|d_{\tilde{\mathcal{G}}}(0,\tilde{X}_{\lfloor tn\rfloor})-
d_{\tilde{\mathcal{G}}}(0,\pi(Z_{\lfloor tn\rfloor}))\right|+
\sup_{t\in[0,T]}\left|d_{\tilde{\mathcal{G}}}(0,\pi(Z_{\lfloor tn\rfloor}))-\delta(d)|Z_{\lfloor tn\rfloor}|\right|\\
&\leq & \sup_{|m|\leq Z_n^*+1} {\rm diam}\tilde{\mathcal{G}}_m+ \sup_{|m|\leq Z_n^*}\left| d_{\tilde{\mathcal{G}}}(0,C_m)-\delta(d)m\right|,
\end{eqnarray*}
where ${\rm diam}\tilde{\mathcal{G}}_m$ is the diameter of the graph $\tilde{\mathcal{G}}_m$, as defined in Section \ref{ergsec}, which is bounded above by $T_{m+1}-T_m$. Now, from the convergence results above, we can conclude that the sequence $(n^{-1/2}\delta(d)Z_n^*)_{n\geq 0}$ converges in distribution under $\tilde{\mathbf{P}}_0^{\tilde{\mathcal{G}}}$ to
$B^*:=\sup_{t\leq T}|B_{t\kappa_1(d)}|$, which is a finite random variable. Consequently, to prove (\ref{close}) it will be enough to establish that
\begin{equation}\label{inc}
n^{-1}{\sup_{|m|\leq n}(T_{m+1}-T_m)}\rightarrow 0,
\end{equation}
\begin{equation}\label{inc2}
n^{-1}\sup_{|m|\leq n}\left|d_{\tilde{\mathcal{G}}}(0,C_m)-\delta(d)m\right|\rightarrow0,
\end{equation}
as $n\rightarrow\infty$, $\mathbf{P}$-a.s. Both of these limits are easily deduced from Lemma \ref{exponent}, which completes the proof of the first part of the theorem.

For the second part of the theorem, we first note that
\begin{eqnarray*}
\lefteqn{\sup_{t\in[0,T]}\left|\tilde{X}_{\lfloor tn\rfloor}-\tilde{S}_{\lfloor\tau(d)Z_{\lfloor tn\rfloor}\rfloor}\right|}\\
&\leq&
\sup_{t\in[0,T]}\left|\tilde{X}_{\lfloor tn\rfloor}-\pi(Z_{\lfloor tn\rfloor})\right|+\sup_{t\in[0,T]}\left|\pi(Z_{\lfloor tn\rfloor})-
\tilde{S}_{\lfloor\tau(d)Z_{\lfloor tn\rfloor}\rfloor}\right|\\
&\leq & \sup_{|m|\leq Z_n^*+1} \sup_{T_{m}\leq l,l'\leq T_{m+1}}\left|\tilde{S}_l-\tilde{S}_{l'}\right|+\sup_{|m|\leq Z_n^*}\left|\pi(m)-\tilde{S}_{\lfloor \tau(d)m\rfloor}\right|.
\end{eqnarray*}
where $(Z_n^*)_{n\geq 0}$ is defined as above. From (\ref{taud}) and the weak convergence of $(n^{-1/4}\tilde{S}_{\lfloor tn^{1/2}\rfloor})_{t\geq0}$ under $\mathbf{P}$,
\[\lim_{n\rightarrow\infty}\tilde{\mathbb{P}}\left(n^{-1/4}\sup_{|m|\leq Z_n^*}\left|\pi(m)-\tilde{S}_{\lfloor \tau(d)m\rfloor}\right|>\varepsilon\right)=0,\]
where we also apply the fact that $(n^{-1/2}Z_n^*)_{n\geq 0}$ converges in distribution under $\tilde{\mathbb{P}}$ to the finite random variable $B^*$. Furthermore, again applying the weak convergence of $(n^{-1/4}\tilde{S}_{\lfloor tn^{1/2}\rfloor})_{t\geq0}$ and $n^{-1/2}Z_n^*$, (\ref{inc}) implies that
\[\lim_{n\rightarrow\infty}\tilde{\mathbb{P}}\left(n^{-1/4}\sup_{|m|\leq Z_n^*+1} \sup_{T_{m}\leq l,l'\leq T_{m+1}}\left|\tilde{S}_l-\tilde{S}_{l'}\right|>\varepsilon\right)=0.\]
Consequently, we have proved that
\begin{equation}\label{close2}
\lim_{n\rightarrow\infty}\tilde{\mathbb{P}}\left(n^{-1/4}\sup_{t\in[0,T]}\left|\tilde{X}_{\lfloor tn\rfloor}-\tilde{S}_{\lfloor \tau(d)Z_{\lfloor tn\rfloor}\rfloor}\right|>\varepsilon\right)=0.
\end{equation}
From the convergence result for $Z$ described at the beginning of the proof and the weak convergence of $(n^{-1/4}\tilde{S}_{\lfloor tn^{1/2}\rfloor})_{t\geq0}$, it is possible to check that the joint law of
\[\left(n^{-1/4}\tilde{S}_{\lfloor tn^{1/2}\rfloor},n^{-1/2}\tau(d)Z_{\lfloor tn\rfloor}\right)_{t\geq0},\]
under $\tilde{\mathbb{P}}$, converges to the joint law of $(W^{(d)}_t,B_{t\kappa_2(d)})_{t\geq0}$. Composing the two processes of the above pair and applying (\ref{close2}), the second conclusion of the theorem follows.
\end{proof}

To adapt the proof of the above theorem to deal with the case of the random walk $X$ on the range of a single random walk ${\mathcal{G}}$, we start by showing how the process $\tilde{X}$ observed on the set \begin{equation}\label{tildevdef}
V(\tilde{\mathcal{G}})_+:=\{S_n:n\geq T_1\}
\end{equation}
satisfies the conclusions of Theorem \ref{d5thm} by applying a time-change argument, and then complete the proof by demonstrating that the time-changed process can be coupled with $X$ in such a way that it is uniformly close on compact time intervals. First, construct an additive functional $A^Z=(A^Z_n)_{n\geq 0}$ related to the process $Z$, defined as at (\ref{zdef}), by setting $A_0^Z=0$ and
\[A_n^Z:=\sum_{m=0}^{n-1}\mathbf{1}_{\{Z_m\geq 0\}},\]
for $n\geq 1$. The analogous functional $A^B=(A^B_t)_{t\geq 0}$ for the Brownian motion $B$ is obtained by setting
\[A^B_t:=\int_0^t\mathbf{1}_{\{B_s\geq 0\}}ds,\]
for $t\geq 0$.

{\lem\label{zlem} Let $d\geq 5$. For $\mathbf{P}$-a.e. realisation of $\tilde{\mathcal{G}}$, the joint law of
\[\left(n^{-1/2}\delta(d)Z_{\lfloor tn\rfloor},n^{-1}A^Z_{\lfloor tn\rfloor}\right)_{t\geq0}\]
under $\mathbf{P}_0^{\tilde{\mathcal{G}}}$ converges to the joint law  of $(B_{t\kappa_1{(d)}}, \kappa_1(d)^{-1}A^B_{t\kappa_1(d)})_{t\geq 0}$.}
\begin{proof} Similarly to the previous proof, the convergence of the first coordinate follows from Corollary \ref{coro}. The lemma is a straightforward consequence of this result and the fact that $\int_0^t\mathbf{1}_{\{B_s\in[-\varepsilon,\varepsilon]\}}ds\rightarrow 0$ almost surely as $\varepsilon\rightarrow 0$, for every finite $t$.
\end{proof}

The time-change of $\tilde{X}$ that we consider will be based on the additive functional $A^{\tilde{X}}=(A^{\tilde{X}}_n)_{n\geq 0}$, defined by setting $A_0^{\tilde{X}}:=0$ and \[A_n^{\tilde{X}}:=\sum_{m=0}^{n-1}\mathbf{1}_{\{\tilde{X}_m,\tilde{X}_{m+1}\in V(\tilde{\mathcal{G}})_+\}},\]
for $n\geq 1$, where $V(\tilde{\mathcal{G}})_+$ was introduced at (\ref{tildevdef}). That $A^{\tilde{X}}$ and $A^Z$ are close is demonstrated by the following lemma.

{\lem \label{axlem}Let $d\geq 5$. For $\mathbf{P}$-a.e. realisation of $\tilde{\mathcal{G}}$, for every $\varepsilon, T\in(0,\infty)$,
\[\lim_{n\rightarrow\infty}\mathbf{P}_0^{\tilde{\mathcal{G}}}\left(n^{-1}\sup_{m\leq Tn}\left|A^{\tilde{X}}_m-A^Z_m\right|>\varepsilon\right)=0.\]}
\begin{proof} By definition, if $n\in[H_m,H_{m+1})$, then $H^{-1}_n=m+1$. It follows that the condition $Z_n\geq 3$, which is equivalent to \[\tilde{X}_{H_{H^{-1}_n}}\in\{S_m:m\geq T_3\},\]
implies that $\tilde{X}_n,\tilde{X}_{n+1}\in V(\tilde{\mathcal{G}})_+$. Conversely, one can check that $\tilde{X}_n,\tilde{X}_{n+1}\in V(\tilde{\mathcal{G}})_+$ implies $Z_n\geq 0$. Consequently,
\[\sup_{m\leq Tn}\left|A^{\tilde{X}}_m-A^Z_m\right|\leq \sum_{m=0}^{n}\mathbf{1}_{\{Z_m\in\{0,1,2\}\}}.\]
By applying the convergence of the rescaled $Z$ to Brownian motion described in Lemma \ref{zlem}, the result is readily deduced from this bound.
\end{proof}

We now introduce a process $\tilde{X}^+=(\tilde{X}^+_n)_{n\geq 0}$ by setting $\tilde{X}^+_n:=\tilde{X}_{\tilde{\alpha}(n)}$, where $\tilde{\alpha}=(\tilde{\alpha}(n))_{n\geq 0}$ is the discrete time inverse of $A^{\tilde{X}}$, defined to satisfy
\[\tilde{\alpha}(n):=\max\{m:A^{\tilde{X}}_m\leq n\}.\]
Note that $\tilde{X}^+$ has the same distribution as the simple random walk on the graph generated by $(S_m)_{m\geq T_1}$ started from $C_1$. We can construct an identically distributed process, $X^+$ say, from $X$. To this end, first let $A^X$ be an additive function defined from $X$, analogously to the definition of $A^{\tilde{X}}$ from $\tilde{X}$. Take $\alpha$ to be the discrete time inverse of $A^X$ and then set $X^+_n:=X_{\alpha(n)}$ for $n\geq 0$. To establish that $X^+$ and $\tilde{X}^+$ have the same distribution is an elementary exercise, and this equivalence allows us to prove the following result.

{\lem\label{prevlem} Theorem \ref{d5thm} holds with $X^+$ in place of $X$.}
\begin{proof} From (\ref{close}) and Lemmas \ref{zlem} and \ref{axlem}, it is possible to deduce that for $\mathbf{P}$-a.e. realisation of $\tilde{\mathcal{G}}$, the joint law of
\[\left(n^{-1/2}d_{\tilde{\mathcal{G}}}(0,\tilde{X}_{\lfloor tn\rfloor}),n^{-1}A^{\tilde{X}}_{\lfloor tn\rfloor}\right)_{t\geq 0}\]
under $\mathbf{P}_0^{\tilde{G}}$ converges to the joint law of $(|B_{t\kappa_1(d)}|,\kappa_1(d)^{-1}A^B_{t\kappa_1(d)})_{t\geq0}$. Since the process $(B_{\alpha^B(t\kappa_1(d))})_{t\geq 0}$ has the same distribution as
$(|B_{t\kappa_1(d)}|)_{t\geq 0}$, where $\alpha^B$ is the right-continuous inverse of $A^B$, it easily follows that the law of
\[\left(n^{-1/2}d_{\tilde{\mathcal{G}}}(0,\tilde{X}^+_{\lfloor tn\rfloor})\right)_{t\geq 0}\]
under $\mathbf{P}_0^{\tilde{G}}$ converges to the law of $(|B_{t\kappa_1(d)}|)_{t\geq0}$. Taking into account the comments preceding this lemma, there is no problem in substituting $X^+$ for $\tilde{X}^+$. To replace $d_{\tilde{\mathcal{G}}}$ by $d_{\mathcal{G}}$, it suffices to note that
\[\sup_{x\in V({\mathcal{G}})}\left|d_{\mathcal{G}}(0,x)-d_{\tilde{\mathcal{G}}}(0,x)\right|\leq T_1-T_0,\]
which is finite, $\mathbf{P}$-a.s. We have therefore established the first convergence result that we are required to prove. By applying (\ref{close2}) in place of (\ref{close}), the second convergence result can be proved similarly.
\end{proof}

To complete the proof of Theorem \ref{d5thm}, all that remains to show is that $X^+$ is a good approximation for the simple random walk $X$. We do this in the next lemma by suitably bounding the amount of time that the sample paths of $X$ spend close to 0. In doing so, we introduce a measure $\mu_\mathcal{G}$ on $V(\mathcal{G})$ by defining $\mu_{\mathcal{G}}(\{x\}):={\rm deg}_{\mathcal{G}}(x)$ for $x\in V(\mathcal{G})$, and consider the ball $B_\mathcal{G}(x,r):=\{y\in V(\mathcal{G}):d_\mathcal{G}(x,y)\leq r\}$.

{\lem\label{xplusapprox} Let $d\geq 5$ and $\varepsilon, T>0$. For $\mathbf{P}$-a.e. realisation of $\mathcal{G}$,
\begin{equation}\label{approx1}
\lim_{n\rightarrow\infty}\mathbf{P}_0^\mathcal{G}\left(n^{-1/2}\sup_{m\leq Tn}\left|d_\mathcal{G}(0,X_m)-d_{\mathcal{G}}(0,X^+_m)\right|>\varepsilon\right)=0.
\end{equation}
Moreover,
\begin{equation}\label{approx2}
\lim_{n\rightarrow\infty}\mathbb{P}\left(n^{-1/4}\sup_{m\leq Tn}\left|X_m-X^+_m\right|>\varepsilon\right)=0.
\end{equation}}
\begin{proof} We clearly have that $\mu_\mathcal{G}(B_\mathcal{G}(S_m,n))\geq n$, for every $m,n\geq 0$. By applying the argument of \cite{BCK}, Proposition 3.3 (and also the Cauchy-Schwarz inequality as in the proof of \cite{BCK}, Proposition 3.8), it follows that for $\mathbf{P}$-a.e. realisation of $\mathcal{G}$, there exists a deterministic finite constant $c_1$ such that
\begin{equation}\label{tpupper}
\mathbf{P}_0^\mathcal{G}(X_n\in\{S_m:0\leq m\leq T_1\})\leq c_1n^{-1/2},
\end{equation}
and, therefore, there exists a finite constant $c_2$ such that
\begin{eqnarray*}
\mathbf{P}_0^\mathcal{G}\left(n^{-1}\sup_{m\leq n}\left|A^X_m-m\right|>\varepsilon\right)&\leq & \varepsilon^{-1}n^{-1}\mathbf{E}_0^\mathcal{G}\sum_{m=0}^{n-1}\mathbf{1}_{\{\{X_m,X_{m+1}\}\cap\{S_l:0\leq l\leq T_1\}\neq \emptyset\}}\\
&\leq & 2\varepsilon^{-1}n^{-1}\sum_{m=0}^{n}\mathbf{P}_0^\mathcal{G}\left(X_m\in\{S_l:0\leq l\leq T_1\}\right)\\
&\leq& c_2n^{-1/2},
\end{eqnarray*}
for $n\geq 1$. This bound implies that, for every $\varepsilon, T\in(0,\infty)$,
\[\lim_{n\rightarrow \infty}\mathbb{P}\left(n^{-1}\sup_{0\leq m\leq Tn}\left|A^X_m-m\right|>\varepsilon\right)=0,\]
and, by the continuous mapping theorem, the same result holds with $\alpha$ in place of $A^X$. Consequently, since $X^+_n=X_{\alpha(n)}$ by definition, to complete the proof it will suffice to establish the tightness of $\{(n^{-1/2}d_{\mathcal{G}}(0,X_{\lfloor tn\rfloor}))_{t\geq 0}\}_{n\geq 1}$ and $\{(n^{-1/4}X_{\lfloor tn\rfloor})_{t\geq 0}\}_{n\geq 1}$ under the appropriate measures.

Suppose that, for some $\delta, \varepsilon>0$ and $n\geq 1$,
\[\sup_{\substack{l,m\leq Tn \\ |l-m|<\delta n}}\left|d_\mathcal{G}(0,X_l)-d_{\mathcal{G}}(0,X_m)\right|>2\varepsilon n^{1/2}\]
and $\sup_{0\leq m\leq T_1}d_\mathcal{G}(0,S_m)<\varepsilon n^{1/2}$. Under these conditions, there exist $l,m\leq Tn$ with $|l-m|<\delta n$ such that $|d_\mathcal{G}(0,X_l)-d_{\mathcal{G}}(0,X_m)|>\varepsilon n^{1/2}$ and also $d_\mathcal{G}(0,X_l),d_\mathcal{G}(0,X_m)>\varepsilon n^{1/2}$. By the definition of $\alpha$, the lower bound on $d_\mathcal{G}(0,X_l)$ and $d_\mathcal{G}(0,X_m)$ implies that $l=\alpha(u), m=\alpha(v)$ for some $u,v\geq 0$, and for this choice of $u,v$ we clearly have that
$|d_\mathcal{G}(0,X^+_u)-d_{\mathcal{G}}(0,X^+_v)|>\varepsilon n^{1/2}$. Now, observing that $|A^X_{m'}-A^X_{n'}|\leq |m'-n'|$ for every $m',n'\geq 0$, so that $|\alpha(m')-\alpha(n')|\geq |m'-n'|$, it is possible to deduce that $u,v\leq Tn$ and $|u-v|\leq|\alpha(u)-\alpha(v)|=|l-m|<\delta n$. Thus, for $\mathbf{P}$-a.e. realisation of $\mathcal{G}$ and every $\delta,\varepsilon, T\in(0,\infty)$,
\begin{eqnarray}
\lefteqn{\mathbf{P}_0^\mathcal{G}\left(\sup_{\substack{l,m\leq Tn \\ |l-m|<\delta n}}\left|d_\mathcal{G}(0,X_l)-d_{\mathcal{G}}(0,X_m)\right|>2\varepsilon n^{1/2}\right)}\nonumber\\
&\leq&\mathbf{P}_0^\mathcal{G}\left(\sup_{\substack{l,m\leq Tn \\ |l-m|<\delta n}}\left|d_\mathcal{G}(0,X^+_l)-d_{\mathcal{G}}(0,X^+_m)\right|>\varepsilon n^{1/2}\right)\label{star}
\end{eqnarray}
for large $n$, where we apply the fact that $\sup_{0\leq m\leq T_1}d_\mathcal{G}(0,S_m)<\varepsilon n^{1/2}$ for large $n$, $\mathbf{P}$-a.s. By Lemma \ref{prevlem}, for $\mathbf{P}$-a.e. realisation of $\mathcal{G}$, the sequence $(n^{-1/2}d_\mathcal{G}(0,X^+_{\lfloor tn\rfloor}))_{t\geq 0}$ is convergent in distribution as $n\rightarrow \infty$, and therefore tight, under $\mathbf{P}_0^\mathcal{G}$. Hence (\ref{star}) implies that $(n^{-1/2}d_\mathcal{G}(0,X_{\lfloor tn\rfloor}))_{t\geq 0}$ is also tight under $\mathbf{P}_0^\mathcal{G}$, for $\mathbf{P}$-a.e. $\mathcal{G}$, which establishes (\ref{approx1}). The proof of (\ref{approx2}) is similar and is omitted.
\end{proof}

Combining the two previous lemmas yields Theorem \ref{d5thm}, and we now prove Theorem \ref{d5hk}. The effective resistance operator on $\mathcal{G}$ will be denoted by $R_\mathcal{G}$.

\begin{proof}[Proof of Theorem \ref{d5hk}] To prove the quenched transition probability asymptotics of (\ref{qtp}), by applying ideas from \cite{BCK} and \cite{Kumagai} (cf. \cite{KM}, Section 3) it will suffice to demonstrate that there exist deterministic constants $c_1,c_2,c_3,c_4\in(0,\infty)$ such that, for $\mathbf{P}$-a.e. realisation of $\mathcal{G}$,
\begin{equation}\label{vol}
c_1n\leq \mu_\mathcal{G}\left(B_\mathcal{G}(0,n)\right)\leq c_2n,
\end{equation}
\begin{equation}\label{res}
c_3n\leq R_\mathcal{G}\left(0,B_\mathcal{G}(0,n)^c\right)\leq c_4n,
\end{equation}
for large $n$, where $\mu_\mathcal{G}$ and $B_\mathcal{G}$ were defined above Lemma \ref{xplusapprox}. The lower bound of (\ref{vol}) is obvious. To prove the upper bound, first observe that $d_\mathcal{G}(0,C_n)\geq n$ for every $n\geq 1$, $\mathbf{P}$-a.s. Consequently, $B_\mathcal{G}(0,n)\subseteq \{S_m:0\leq m\leq T_n\}$. Furthermore, it is straightforward to check that $\mu_\mathcal{G}(\{S_m:0\leq m\leq T_n\})\leq 2T_n+2$. Thus $\mu_\mathcal{G}(B_\mathcal{G}(0,n))\leq 4T_n$, and the upper bound at (\ref{vol}) follows from the ergodic limit result for the cut-times $T_n$ at (\ref{taud}). By the connectedness of the graph $\mathcal{G}$, there exists at least one path of length $n+1$ from 0 to $B_\mathcal{G}(0,n)^c$. This readily implies that $R_\mathcal{G}(0,B_\mathcal{G}(0,n)^c)\leq n+1\leq 2n$, for every $n$. To complete the proof of (\ref{res}), we first observe that
\[\sup_{m\leq T_n}d_\mathcal{G}(0,S_m)\leq \sup_{0\leq m\leq n}d_{\mathcal{G}}(0,C_m)+\sup_{0\leq m\leq n}(T_{m+1}-T_m)\leq 2\delta(d)n\]
for large $n$, $\mathbf{P}$-a.s., which can be proved by applying (\ref{inc}) and (\ref{inc2}).
Note that to apply (\ref{inc2}), we need to replace $d_{\tilde{\mathcal{G}}}$ by $d_\mathcal{G}$, which can be justified as in the proof of Lemma \ref{prevlem}. Thus $\{S_m:0\leq m\leq T_n\}\subseteq B_\mathcal{G}(0,2\delta(d)n))$ for large $n$, $\mathbf{P}$-a.s. This result implies that, $\mathbf{P}$-a.s., for large $n$, any path from 0 to $B_\mathcal{G}(0,2\delta(d)n))^c$ must pass through the edges $\{S_{T_m},S_{T_m+1}\}$, $1\leq m\leq n-1$. Consequently, simple properties of resistance yield that $R_\mathcal{G}(0,B_\mathcal{G}(0,2\delta(d)n))^c)\geq n-1$ for large $n$, $\mathbf{P}$-a.s., and the lower bound at (\ref{res}) follows.
\end{proof}

We conclude this section with a couple of further remarks about the results in high dimensions. Firstly, in addition to the quenched transition probability bounds of Theorem \ref{d5hk}, it is also possible to deduce corresponding annealed bounds. More specifically, there exist deterministic constants $c_1,c_2\in(0,\infty)$ such that
\begin{equation}\label{anntp}
c_1n^{-1/2}\leq \mathbb{P}(X_{2n}=0)\leq c_2n^{-1/2},
\end{equation}
for every $n\in\mathbb{N}$. For the upper bound we can simply take expectations in (\ref{tpupper}). For the lower bound we apply Fatou's lemma and the lower bound at (\ref{qtp}) to deduce that
\[\liminf_{n\rightarrow\infty}n^{1/2}\mathbb{P}(X_{2n}=0)\geq \mathbf{E}\left(\liminf_{n\rightarrow\infty}n^{1/2}\mathbf{P}_0^\mathcal{G}(X_{2n}=0)\right)\geq c_3,\]
where $c_3$ is a strictly positive deterministic constant, and the result follows.

Secondly, for the alternative version of the random walk $X$ described in the introduction with transition probabilities $P_\mathcal{G}(x,y):=\mu_{xy}/\mu_x$, where $\mu_{xy}$ is defined at (\ref{muxy}) and $\mu_x:=\sum_{y:\{x,y\}\in E(\mathcal{G})}\mu_{xy}$, essentially the same proofs will yield the results corresponding to Theorem \ref{d5thm} and \ref{d5hk}. However in the electrical network interpretation of the random walk on $\tilde{\mathcal{G}}$, the resistance metric $R_{\tilde{\mathcal{G}}}$ we need to consider is the effective resistance resulting when each pair of adjacent vertices $x,y\in V(\tilde{\mathcal{G}})$ is connected by a wire of conductance
\[\tilde{\mu}_{xy}:=\#\left\{n:\{\tilde{S}_n,\tilde{S}_{n+1}\}=\{x,y\}\right\}.\]
We also need to replace ${\rm deg}_{\tilde{\mathcal{G}}}(x)$ by $\tilde{\mu}_x:=\sum_{y:\{x,y\}\in E(\tilde{\mathcal{G}})}\tilde{\mu}_{xy}$ in the definitions of $\mu^{\tilde{\mathcal{G}}}$, $\nu(d)$ and $\eta(d)$. The one point that requires a little more careful checking is that the revised expression for $\eta(d)$,
\begin{equation}\label{reviseeta}
\frac{\hat{\mathbf{E}}\left(\tilde{\mu}_0
\mathbf{E}^{\tilde{\mathcal{G}}}_0 H_1\right)}{\hat{\mathbf{E}}\tilde{\mu}_0},
\end{equation}
is finite. To do this, we proceed similarly to (\ref{bound1}) to deduce that
\begin{eqnarray*}
\tilde{\mu}_0\mathbf{E}_0^{\tilde{\mathcal{G}}}H_1&\leq &\tilde{\mu}_0+
\sum_{\substack{x\in V(\tilde{\mathcal{G}})\backslash \{C_{-1},C_1\}:\\\{0,x\}\in E(\tilde{\mathcal{G}})}}\tilde{\mu}_{0x}\mathbf{E}_x^{\tilde{\mathcal{G}}}H_1\\
&\leq &\tilde{\mu}_0+\sum_{\substack{x\in V(\tilde{\mathcal{G}})\backslash \{C_{-1},C_1\}:\\\{0,x\}\in E(\tilde{\mathcal{G}})}}\tilde{\mu}_{0x}R_{\tilde{\mathcal{G}}}(0,x)\mu_{\tilde{\mathcal{G}}}\left(A\backslash\{0\}\right)\\
&\leq &\tilde{\mu}_0+2d\mu_{\tilde{\mathcal{G}}}\left(A\backslash\{0\}\right)\\
&\leq &(1+4d)(T_{1}-T_{-1}),
\end{eqnarray*}
where $A$ is defined as in the proof of Lemma \ref{hitting} and for the third inequality we apply the fact that $R_{\tilde{\mathcal{G}}}(0,x)\leq \tilde{\mu}_{0x}^{-1}$. The finiteness of (\ref{reviseeta}) follows.

\section{Behaviour at the critical dimension $d=4$}\label{d4sec}

In this section we will prove Theorems \ref{d4thm} and \ref{tdthm}, demonstrating that when $d=4$ the process $X$ and its transition density do not satisfy the same scaling results as in the high-dimensional case, exhibiting logarithmic corrections to the leading order polynomial behaviour. We start by stating some known properties for the random walk $S$ that will be used to establish properties of the range of the random walk $\mathcal{G}$. In a change of notation from the previous section we will write $(T_n)_{n\geq 1}$ to represent the elements of the set of cut-times $\mathcal{T}$, as defined at (\ref{osct}), arranged in an increasing order. Denoting the loop-erasure of $(S_m)_{m=0}^n$ by $L^n$ (see \cite{Lawlerbook}, Section 7.2, for a definition), we define $Y=(Y_n)_{n\geq0}$ by setting $Y_n$ to be equal to the number of edges in the path $L^n$ (so that in the notation of \cite{Lawlerbook}, $Y_n=\rho_n(n)$).

{\lem \label{rwprops}Let $d=4$.\\
(a) There exists a deterministic constant $c_1\in(0,\infty)$ such that
\[\frac{T_n}{n(\ln n)^{1/2}}\buildrel \mathbf{P}\over\rightarrow c_1,\]
as $n\rightarrow\infty$.\\
(b) The process $Y$ satisfies
\[\lim_{\lambda\rightarrow\infty}\limsup_{n\rightarrow\infty}\mathbf{P}\left( \sup_{m\leq n} Y_m\geq \lambda n (\ln n)^{-1/3}\right)=0.\]
(c) $\mathbf{P}$-a.s., the simple random walk $S$ satisfies
\[\frac{\#\{S_m:0\leq m\leq n\}}{n}\rightarrow1,\]
as $n\rightarrow\infty$.
}
\begin{proof} Part (a) follows from the discussion in the introduction of \cite{LawlerEJP}. To prove part (b), we consider the process $Y'=(Y'_n)_{n\geq 0}$, defined by letting $Y_n'$ be equal to the number of the first $n$ points of $S$ retained after loop-erasing the whole path $(S_{m})_{m\geq 0}$ (in \cite{Lawlerbook}, this is the process $\rho(n)$). The asymptotic behaviour of the expectation of $Y'_n$ is given by
\[\mathbf{E}Y'_n\sim c_2n(\ln n)^{-1/3},\]
as $n\rightarrow\infty$, for some deterministic constant $c_2\in(0,\infty)$ (see \cite{LERW}, for example). Furthermore, it is possible to deduce that
\[\sup_{t\in[0,1]}\left|\frac{Y'_{\lfloor tn\rfloor}}{\mathbf{E}Y'_n}-t\right|\buildrel \mathbf{P}\over\rightarrow 0,\]
as $n\rightarrow\infty$, by making a simple adaptation to the analogous result for the ``inverse'' of $Y'$ proved in \cite{Lawlerbook}, Section 7.7. Consequently, part (b) of the lemma will follow if we can establish that
\begin{equation}\label{needed}
\frac{\sup_{m\leq n}|Y_m-Y'_m|}{n(\ln n)^{-1/3}}\buildrel \mathbf{P}\over\rightarrow 0,
\end{equation}
as $n\rightarrow\infty$. First, for $n$ large enough, choose $0=j_0<j_1<\dots<j_k=n$ such that
\[\tfrac{1}{2}n(\ln n)^{-2}\leq j_i-j_{i-1} \leq 2n(\ln n)^{-2},\]
and $k\sim (\ln n)^2$, and define for $i=1,\dots k$,
\[Z_{i}=\mathbf{1}\{\mathcal{T}\cap[j_i-j_i(\ln j_i)^{-6},j_i]=\emptyset\},\]
where $\mathcal{T}$ is the set of cut-times of $S$. As is demonstrated by \cite{Lawlerbook}, Lemma 7.7.4, we can check that there exist constants such that (for large enough $n$)
\[\mathbf{E}Z_i\leq c_3\frac{\ln\ln j_i}{\ln j_i}\leq c_4\frac{\ln\ln n}{\ln n},\hspace{20pt}i=1,2,\dots,k.\]
By simple considerations of the structure of the path, we can check that
\[\sup_{m\leq n}|Y_{m}-Y'_m|\leq n(\ln n)^{-6} + 2n(\ln n)^{-2}\left(1+\sum_{i=1}^k Z_i\right).\]
Thus, for $\varepsilon>0$,
\begin{eqnarray*}
\lefteqn{\limsup_{n\rightarrow\infty}\mathbf{P}\left(\sup_{m\leq n}|Y_m-Y'_m|\geq \varepsilon n(\ln n)^{-1/3}\right)}\\
&\leq &\limsup_{n\rightarrow\infty}\mathbf{P}\left(4n(\ln n)^{-2}\sum_{i=1}^k Z_i\geq \varepsilon n(\ln n)^{-1/3}\right)\\
&\leq &\limsup_{n\rightarrow\infty}\frac{4\sum_{i=1}^k\mathbf{E}Z_i}{(\ln n)^{5/3}}\\
&\leq &\limsup_{n\rightarrow\infty}\frac{c \ln\ln n}{(\ln n)^{2/3}}\\
&=&0,
\end{eqnarray*}
which confirms (\ref{needed}). Finally, part (c) is an easy consequence of \cite{ET2}, Theorem 12.
\end{proof}

These properties allow us to establish bounds for the volume of a ball centered at 0 with respect to the measure $\mu_\mathcal{G}$, which was defined above Lemma \ref{xplusapprox}. In the proof, we apply the notation $C_n:=S_{T_n}$.

{\lem\label{d4vol} If $d=4$, then
\begin{equation}\label{f}
\lim_{\lambda\rightarrow\infty}\liminf_{n\rightarrow\infty}\mathbf{P}\left(\lambda^{-1} n(\ln n)^{1/3}\leq \mu_\mathcal{G}(B_{\mathcal{G}}(0,n))\leq \lambda n(\ln n)^{1/2}\right)=1.
\end{equation}}
\begin{proof} Since $d_\mathcal{G}(0,S_m)>d_\mathcal{G}(0,C_n)\geq n-1$ for every $m>T_n$, $n\geq 1$, we have that
\begin{equation}\label{inclusion}
B_\mathcal{G}(0,n-1)\subseteq\{S_m:0\leq m\leq T_n\},
\end{equation}
for every $n\geq 1$. Therefore
\begin{eqnarray*}
\lefteqn{\lim_{\lambda\rightarrow\infty}\limsup_{n\rightarrow\infty}\mathbf{P}\left(\mu_\mathcal{G}(B_{\mathcal{G}}(0,n))\geq \lambda n(\ln n)^{1/2}\right)}\\
&\leq&\lim_{\lambda\rightarrow\infty}\limsup_{n\rightarrow\infty}\mathbf{P}\left(\mu_\mathcal{G}(\{S_m:0\leq m\leq T_{n+1}\})\geq\lambda n(\ln n)^{1/2}\right)\\
&\leq &\lim_{\lambda\rightarrow\infty}\limsup_{n\rightarrow\infty}\mathbf{P}\left(4T_{n+1}
\geq\lambda n(\ln n)^{1/2}\right)\\
&=&0,
\end{eqnarray*}
where we apply the observation that $\mu_\mathcal{G}(\{S_m:0\leq m\leq T_n\})\leq 2T_n+2$, which was also used in the proof of Theorem \ref{d5hk}, to deduce the second inequality and Lemma \ref{rwprops}(a) to deduce the final limit. This completes the proof of the right-hand inequality of (\ref{f}).

To deduce the left-hand inequality of (\ref{f}), first observe that $d_\mathcal{G}(0,S_m)\leq Y_m$ for every $n\geq 1$. Thus, by applying Lemma \ref{rwprops}(b), we obtain that
\begin{eqnarray}
\lefteqn{\lim_{\lambda\rightarrow\infty}\limsup_{n\rightarrow\infty}\mathbf{P}\left( \sup_{m\leq n} d_\mathcal{G}(0,S_m)\geq \lambda n (\ln n)^{-1/3}\right)}\nonumber\\
&\leq &\lim_{\lambda\rightarrow\infty}\limsup_{n\rightarrow\infty}\mathbf{P}\left( \sup_{m\leq n} Y_m\geq \lambda n (\ln n)^{-1/3}\right)\nonumber\\
&=&0.\label{dglim}
\end{eqnarray}
Consequently,
\begin{eqnarray*}
\lefteqn{\lim_{\lambda\rightarrow\infty}\limsup_{n\rightarrow\infty}\mathbf{P}\left( \mu_\mathcal{G}(B_\mathcal{G}(0,\lambda n(\ln n)^{-1/3}))\leq \lambda^{-1} n\right)}\\
&\leq & \lim_{\lambda\rightarrow\infty}\limsup_{n\rightarrow\infty}\mathbf{P}\left( \#\{S_m:0\leq m\leq n\}\leq \lambda^{-1} n\right)\\
&=&0,
\end{eqnarray*}
where we apply the fact that $\mu_\mathcal{G}(\{S_m:0\leq m\leq n\})\geq \#\{S_m:0\leq m\leq n\}$ to deduce the inequality and note that Lemma \ref{rwprops}(c) implies the equality. A simple reparameterisation of $n$ and $\lambda$ completes the proof.
\end{proof}

Combining the above results with standard arguments for random walks on graphs allows us to deduce bounds for
\[\tau_\mathcal{G}(0,n):=\inf\{m:X_m\not\in B_\mathcal{G}(0,n)\},\]
the exit time of $X$ from a ball, and its expectation under $\mathbf{P}_0^\mathcal{G}$.

{\lem\label{d4exit} If $d=4$, then
\begin{equation}\label{exit}
\lim_{\lambda\rightarrow\infty}\liminf_{n\rightarrow\infty}\mathbf{P}\left(\lambda^{-1}n^2 \leq \mathbf{E}_0^\mathcal{G}\tau_\mathcal{G}(0,n)\leq \lambda n^2(\ln n)^{1/2}\right)=1,
\end{equation}
and also
\begin{equation}\label{exitnon}
\lim_{\lambda\rightarrow\infty}\liminf_{n\rightarrow\infty}\mathbb{P}\left(\lambda^{-1} n^2(\ln n)^{-4/3}\leq \tau_\mathcal{G}(0,n)\leq \lambda n^2(\ln n)^{1/2}\right)=1.
\end{equation}}
\begin{proof} That $\mathbf{E}_0^\mathcal{G}\tau_\mathcal{G}(0,n)\leq R_\mathcal{G}(0,B_\mathcal{G}(0,n)^c)\mu_\mathcal{G}(B_\mathcal{G}(0,n))$ is a well-known result, see \cite{Telcs}, Lemma 3.6, for example. As in the proof of Theorem \ref{d5hk}, we have that $R_{\mathcal{G}}(0,B_\mathcal{G}(0,n)^c)\leq n+1$, and therefore the right-hand inequality of (\ref{exit}) is a straightforward consequence of Lemma \ref{d4vol}. To prove the left-hand inequality of (\ref{exit}), first suppose that
\begin{equation}\label{supd}
\sup_{m\leq T_{2n}}d_\mathcal{G}(0,S_m)\leq \lambda n (\ln n)^{1/6},
\end{equation}
\begin{equation}\label{2}
\mu_\mathcal{G}(B_\mathcal{G}(0,n))\geq \lambda^{-1} n (\ln n)^{1/3},
\end{equation}
and let $g_B$ be the quenched occupation density of the random walk on $\mathcal{G}$ killed on exiting $B:=B_\mathcal{G}(0,\lambda n(\ln n)^{1/6})$, so that
\[\mathbf{E}_0^\mathcal{G}\tau_\mathcal{G}(0,\lambda n(\ln n)^{1/6})=\sum_{x\in\mathcal{G}}g_B(x)\mu_\mathcal{G}(\{x\}).\]
By (\ref{supd}), any path from 0 to $B^c$ passes through $C_{2n}$, and therefore $R_\mathcal{G}(0,B^c)\geq R(0,C_{2n})\geq 2n-1$. Thus, applying an argument from the proof of \cite{BCK}, Proposition 3.4, for example, it is possible to deduce that $g_B(x)\geq cn$, for every $x\in B_\mathcal{G}(0,n)$ and $n\geq 2$, where $c$ is a strictly positive deterministic constant. Consequently if (\ref{2}) also holds, then
\[\mathbf{E}_0^\mathcal{G}\tau_\mathcal{G}(0,\lambda n(\ln n)^{1/6})\geq c \lambda^{-1}n^2(\ln n)^{1/3},\]
and so, in view of Lemma \ref{d4vol}, to complete the proof of (\ref{exit}) it remains to establish that
\begin{equation}\label{supdlim}
\lim_{\lambda\rightarrow\infty} \limsup_{n\rightarrow\infty}\mathbf{P}\left( \sup_{m\leq T_{2n}}d_\mathcal{G}(0,S_m)\geq \lambda n (\ln n)^{1/6}\right)=0.
\end{equation}
Applying Lemma \ref{rwprops}(a), this is a straightforward adaptation of the result proved at (\ref{dglim}).

On noting that, for $\lambda>0$, $\varepsilon\in(0,1)$,
\begin{eqnarray*}
\mathbb{P}\left(\tau_{\mathcal{G}}(0,n)\geq \lambda\right)&\leq&\mathbf{E}\left(\frac{\mathbf{E}_0^\mathcal{G}\tau_\mathcal{G}(0,n)}{\lambda}\wedge 1\right)\\
&\leq&\mathbf{E}\left(\mathbf{1}_{\{\mathbf{E}_0^\mathcal{G}\tau_\mathcal{G}(0,n)\geq \varepsilon\lambda\}}\right)+\mathbf{E}\left(\frac{\mathbf{E}_0^\mathcal{G}\tau_\mathcal{G}(0,n)}{\lambda}\mathbf{1}_{\{\mathbf{E}_0^\mathcal{G}\tau_\mathcal{G}(0,n)< \varepsilon\lambda\}}\right)\\
&\leq&\mathbf{P}\left(\mathbf{E}_0^\mathcal{G}\tau_\mathcal{G}(0,n)\geq \varepsilon\lambda\right)+\varepsilon,
\end{eqnarray*}
the right-hand inequality of (\ref{exitnon}) is readily deduced from the right-hand inequality of (\ref{exit}). For the left-hand inequality, we adapt the argument of \cite{KM}, Proposition 3.5(a). Firstly, suppose that for some $n,\lambda\geq 2$ and $\varepsilon\in(0,1)$, we have that
\begin{equation}\label{taubounds}
\lambda^{-1}(n\delta)^2\leq \mathbf{E}_0^\mathcal{G}\tau_\mathcal{G}(0,n\delta)\leq\sup_{x}\mathbf{E}_x^\mathcal{G}\tau_\mathcal{G}(0,n\delta)\leq\lambda (n\delta)^2(\ln n\delta)^{1/2},
\end{equation}
where $\delta:=(2\lambda\varepsilon)^{1/2}$, and (\ref{supd}) holds. We also assume that $\delta\in (0,1)$ and $n\delta +1\leq n$. Applying the strong Markov property at $\tau_\mathcal{G}(0,n\delta)$ yields
\begin{eqnarray*}
\lefteqn{\mathbf{P}_0^\mathcal{G}\left(\tau_\mathcal{G}(0,\lambda n (\ln n)^{1/6})\leq \varepsilon n^2\right)}\\
&\leq& \mathbf{P}_0^\mathcal{G}\left(\tau_\mathcal{G}(0,n\delta)< \varepsilon n^2\right)\sup_{x\in B_\mathcal{G}(0,n\delta+1)}\mathbf{P}_x^\mathcal{G}\left(\tau_\mathcal{G}(0,\lambda n (\ln n)^{1/6})\leq \varepsilon n^2\right).
\end{eqnarray*}
To bound the second term, we again apply the strong Markov property, this time at $T_0:=\inf\{m:X_m=0\}$, to obtain that \begin{eqnarray*}
\lefteqn{\mathbf{P}_x^\mathcal{G}\left(\tau_\mathcal{G}(0,\lambda n (\ln n)^{1/6})\leq\varepsilon n^2\right)}\\
&\leq &\mathbf{P}_x^\mathcal{G}\left(\tau_\mathcal{G}(0,\lambda n (\ln n)^{1/6})<T_0\right)+\mathbf{P}_0^\mathcal{G}\left(\tau_\mathcal{G}(0,\lambda n (\ln n)^{1/6})\leq \varepsilon n^2\right),
\end{eqnarray*}
and therefore
\begin{eqnarray}
\lefteqn{\mathbf{P}_0^\mathcal{G}\left(\tau_\mathcal{G}(0,\lambda n (\ln n)^{1/6})\leq \varepsilon n^2\right)}\nonumber\\
&\leq&
 \mathbf{P}_0^\mathcal{G}\left(\tau_\mathcal{G}(0,n\delta)\geq \varepsilon n^2\right)^{-1}
 \sup_{x\in B_\mathcal{G}(0,n\delta+1)}\mathbf{P}_x^\mathcal{G}\left(\tau_\mathcal{G}(0,\lambda n (\ln n)^{1/6})<T_0\right).\label{bound12}
\end{eqnarray}
We now explain how to bound each of these terms. For the first term, by a standard Markov property argument, we have
\begin{equation}\label{mp}
\mathbf{E}_0^\mathcal{G}\tau_\mathcal{G}(0,n\delta)\leq
\varepsilon n^2+\mathbf{P}_0^\mathcal{G}\left(\tau_\mathcal{G}(0,n\delta)\geq\varepsilon n^2\right)\sup_{x\in B_\mathcal{G}(0,n\delta)}\mathbf{E}_x^\mathcal{G}\tau_\mathcal{G}(0,n\delta)
\end{equation}
from which it follows that
\begin{equation}\label{tauprob}
\mathbf{P}_0^\mathcal{G}\left(\tau_\mathcal{G}(0,n\delta)\geq\varepsilon n^2\right)\geq \left( 2\lambda^2(\ln n\delta)^{1/2}\right)^{-1},
\end{equation}
where we apply (\ref{taubounds}) to bound the expectations in (\ref{mp}). The second term of (\ref{bound12}) can be bounded above by
\[\sup_{x\in B_\mathcal{G}(0,n\delta+1)}\frac{R_\mathcal{G}(0,x)}{R_\mathcal{G}(x,B_\mathcal{G}(0,\lambda n (\ln n)^{1/6})^c)}\leq \frac{n\delta+1}{n}=\delta+n^{-1},\]
where the bound on the left-hand side of the above expression is well-known for graphs, see \cite{BGP}, equation (4), for example. To deduce the first inequality, we apply that $R_\mathcal{G}(0,x)\leq d_\mathcal{G}(0,x)\leq n\delta +1$ and also $R_\mathcal{G}(x, B_\mathcal{G}(0,\lambda n (\ln n)^{1/6})^c)\geq R_\mathcal{G}(C_n,C_{2n})\geq n$, where this bound is a consequence of (\ref{supd}), similarly to the resistance bound applied in the proof of the lower bound for the expectation of the exit time from a ball. Combining these bounds, we have proved that $\mathbf{P}_0^\mathcal{G}(\tau_\mathcal{G}(0,\lambda n (\ln n)^{1/6})\leq \varepsilon n^2)\leq \left( 2\lambda^2(\ln n\delta)^{1/2}\right)(\delta+n^{-1})$. Finally, it is possible to deduce from this fact, the left-hand inequality of (\ref{exit}), (\ref{supdlim}) and
\[\lim_{\lambda\rightarrow\infty}\limsup_{n\rightarrow\infty}\mathbf{P}\left(
\sup_x \mathbf{E}_x^\mathcal{G}\tau_\mathcal{G}(0,n)\geq \lambda n^2(\ln n)^{1/2}\right)=0,\]
which can be proved by a simple extension of the proof of the right-hand inequality of (\ref{exit}), that
\[\lim_{\lambda\rightarrow\infty}\limsup_{\varepsilon\rightarrow 0}\limsup_{n\rightarrow\infty}\mathbb{P}\left(\tau_\mathcal{G}(0,\lambda n (\ln n)^{1/6})\leq \varepsilon n^2(\ln n)^{-1}\right)=0,\]
and the result easily follows.
\end{proof}

With these preparations in place, it is now relatively straightforward to prove Theorems \ref{d4thm} and \ref{tdthm}, demonstrating the necessity of logarithmic corrections when $d=4$.

\begin{proof}[Proof of Theorem \ref{d4thm}] For $n\geq 1$, $\lambda>0$, define the events
\begin{eqnarray*}
A_0&:=&\left\{\tau_\mathcal{G}(0,n)\leq \lambda n^2(\ln n)^{1/2}\right\},\\
A_1&:=&\left\{\{S_m:0\leq m\leq  T_{\lfloor\lambda^{-1}n(\ln n)^{-1/6}\rfloor}\}\subseteq B_\mathcal{G}(0,n)\right\},\\
A_2&:=&\left\{\left|C_{\lfloor\lambda^{-1}n(\ln n)^{-1/6}\rfloor}\right|\geq  \lambda^{-2}n^{1/2}(\ln n)^{1/6}\right\}.
\end{eqnarray*}
On the set $A_0\cap A_1\cap A_2$,
\[\max_{m\leq \lambda n^2 (\ln n)^{1/2}}\left|X_m\right|\geq \max_{m\leq \tau_\mathcal{G}(0,n)}\left|X_m\right|\geq \left|C_{\lfloor\lambda^{-1}n(\ln n)^{-1/6}\rfloor}\right|\geq  \lambda^{-2}n^{1/2}(\ln n)^{1/6},\]
where to deduce the second inequality, we note that the definition of a cut-time implies that, on the set $A_1$, by the hitting time $\tau_{\mathcal{G}}(0,n)$, the process $X$ must have hit the vertex $\left|C_{\lfloor\lambda^{-1}n(\ln n)^{-1/6}\rfloor}\right|$. The left-hand inequality of (\ref{d4lim}) will follow easily from this if we can establish that
\[\lim_{\lambda\rightarrow\infty}\limsup_{n\rightarrow\infty}\mathbb{P}\left(A_i^c\right)=0,\]
for $i=0,1,2$. The result for $i=0$ was proved in Lemma \ref{d4exit}. For $i=1$, observe that
\[\mathbf{P}(A_1^c)\leq \mathbf{P}\left(T_{\lfloor\lambda^{-1}n(\ln n)^{-1/6}\rfloor}\geq \lambda^{-1/2}n(\ln n)^{1/3}\right)+\mathbf{P}\left(\max_{m\leq \lambda^{-1/2}n(\ln n)^{1/3}}d_\mathcal{G}(0,S_m)\geq n\right),\]
and so $\lim_{\lambda\rightarrow\infty}\limsup_{n\rightarrow\infty}\mathbf{P}(A_1^c)=0$ by Lemma \ref{rwprops}(a) and (\ref{dglim}). Similarly, $\mathbf{P}(A_2^c)$ is bounded above by
\[\mathbf{P}\left(T_{\lfloor\lambda^{-1}n(\ln n)^{-1/6}\rfloor}\not\in[ \lambda^{-2}n(\ln n)^{1/3},\lambda^{-1/2}n(\ln n)^{1/3}]\right)\]
\[+\mathbf{P}\left(\inf_{m\in[ \lambda^{-2}n(\ln n)^{1/3},\lambda^{-1/2}n(\ln n)^{1/3}]}|S_m|\leq \lambda^{-2}n^{1/2}(\ln n)^{1/6}\right),\]
Again applying Lemma \ref{rwprops}(a) and well-known scaling properties of simple random walk and Brownian motion, we obtain that
\[\lim_{\lambda\rightarrow\infty}\limsup_{n\rightarrow\infty}\mathbf{P}(A_2^c)\leq \lim_{\lambda\rightarrow\infty}\mathbf{P}\left(\inf_{t\in[ \lambda^{-2},\lambda^{-1/2}]}|W^{(d)}_t|\leq \lambda^{-2}\right)\leq \lim_{\lambda\rightarrow\infty}\mathbf{P}\left(\inf_{t\geq 1}|W^{(d)}_t|\leq \lambda^{-1}\right).\]
Since $\inf_{t\geq 1}|W^{(d)}_t|>0$, $\mathbf{P}$-a.s., the right-hand side is equal to 0, as desired.

Similarly, for $n\geq 1$, $\lambda>0$, define
\begin{eqnarray*}
B_0&:=&\left\{\tau_\mathcal{G}(0,n)\geq \lambda^{-1} n^2(\ln n)^{-4/3}\right\},\\
B_1&:=&\left\{T_{n+2}\leq \lambda n (\ln n)^{1/2}\right\},\\
B_2&:=&\left\{\sup_{m\leq \lambda n(\ln n)^{1/2}}|S_m|\leq   \lambda n^{1/2}(\ln n)^{1/4}\right\},
\end{eqnarray*}
so that on the set $B_0\cap B_1\cap B_2$ we have
\begin{eqnarray*}
\max_{m\leq \lambda^{-1} n^2 (\ln n)^{-4/3}}\left|X_m\right|&\leq& \max_{m\leq \tau_\mathcal{G}(0,n)}\left|X_m\right|\leq \sup_{x\in B_\mathcal{G}(0,n+1)}|x|\\
&\leq&
\sup_{m\leq T_{n+2}}|S_m|
\leq \sup_{m\leq \lambda n(\ln n)^{1/2}}|S_m|\leq   \lambda n^{1/2}(\ln n)^{1/4},
\end{eqnarray*}
where we apply (\ref{inclusion}) to deduce the third inequality. Thus the right-hand inequality of (\ref{d4lim}) is a consequence of the fact that
\[\lim_{\lambda\rightarrow\infty}\limsup_{n\rightarrow\infty}\mathbb{P}\left(B_i^c\right)=0,\]
for $i=0,1,2$, which can be deduced by applying Lemma \ref{d4exit} in the case $i=0$, Lemma \ref{rwprops}(a) for $i=1$ and simple random walk scaling properties for $i=2$.
\end{proof}

Before proceeding with our final proof, let us remark that the result of Theorem \ref{d4thm} shows that the extra intersections of $S$ in $d=4$ lead to the random walk $X$ moving logarithmically more quickly away from 0 with respect to Euclidean distance than in higher dimensions. With respect to the graph distance, however, (\ref{exit}) suggests that with respect to the graph distance $X$ will move no quicker when $d=4$ than in higher dimensions. Although this may at first seem paradoxical, it can be explained by observing that the extra connectivity of the graph in $d=4$ allows $X$ to access more easily points that are later in time on the $S$ path, but reduces the graph distance to them.

\begin{proof}[Proof of Theorem \ref{tdthm}] A standard argument, see \cite{KM}, Proposition 3.1(a), for example, implies that if $\mu_\mathcal{G}(B_\mathcal{G}(0,R))\geq \lambda^{-1}R(\ln R)^{1/3}$, then $\mathbf{P}^\mathcal{G}_0(X_{2n}=0)\leq c\lambda n^{1/2}(\ln n)^{-1/6}$ for $\tfrac{1}{2}R^2(\ln R)^{1/3}\leq n\leq R^2(\ln R)^{1/3}$, where $c$ is a finite deterministic constant. Thus the upper transition probability bound of Theorem \ref{tdthm} can be deduced from the lower volume bound of Lemma \ref{d4vol}.

For the lower transition probability bound, another standard argument can be applied. First, if $2R^2\leq \mathbf{E}_0^\mathcal{G}\tau_\mathcal{G}(0,\lambda R)$ and $\mu_\mathcal{G}(B_\mathcal{G}(0,\lambda R))\leq \lambda^2 R(\ln R)^{1/2}$, then, similarly to (\ref{tauprob}), we have that $\mathbf{P}_0^\mathcal{G}\left(\tau_\mathcal{G}(0,\lambda R)> R^2\right)\geq \left( \lambda^3(\ln R)^{1/2}\right)^{-1}$. Applying Cauchy-Schwarz, as in the proof of \cite{KM}, Proposition 3.2, for example, it follows that
\[\mu_\mathcal{G}(B_\mathcal{G}(0,\lambda R)) \mathbf{P}_0^\mathcal{G}(X_{2R^2}=0)\geq \left( \lambda^6(\ln R)\right)^{-1},\]
and, therefore, our upper volume bound assumption implies that \[\mathbf{P}_0^\mathcal{G}(X_{2R^2}=0)\geq \left( \lambda^8R(\ln R)^{3/2}\right)^{-1}.\]
Thus the desired result is a consequence of Lemmas \ref{d4vol} and \ref{d4exit}.
\end{proof}
\hspace{10pt}\\
\textbf{Acknowledgements} The author would like to thank Daisuke Shiraishi for carefully reading an earlier version of this article and pointing out a number of errors that appeared there, and also an anonymous referee for suggesting some additional references.

%\bibliography{david}
%\bibliographystyle{amsplain}

\def\cprime{$'$}
\providecommand{\bysame}{\leavevmode\hbox to3em{\hrulefill}\thinspace}
\providecommand{\MR}{\relax\ifhmode\unskip\space\fi MR }
\providecommand{\MRhref}[2]{%
  \href{http://www.ams.org/mathscinet-getitem?mr=#1}{#2}
}
\providecommand{\href}[2]{#2}

\end{document}